Some thoughts on Le Cam's statistical decision theory

by

David Pollard
Statistics Department
Yale University
Box 208290 Yale Station
New Haven, CT 06520
Email: david.pollard@yale.edu

**Abstract.** *The paper contains some musings about the abstractions introduced by Lucien Le Cam into the asymptotic theory of statistical inference and decision theory. A short, self-contained proof of a key result (existence of randomizations via convergence in distribution of likelihood ratios), and an outline of a proof of a local asymptotic minimax theorem, are presented as an illustration of how Le Cam's approach leads to conceptual simplifications of asymptotic theory.*

## 1. Convergence of experiments

Over a period of almost half a century, Lucien Le Cam developed a general theory for handling asymptotic problems in statistical decision theory. At the core of his thinking was the concept of a distance between statistical models (or experiments, in Le Cam's terminology). When two models are close in Le Cam's sense, there is an automatic transfer of solutions to certain types of decision theoretic problems from one model to the other.

Most of the theory was described, in very general form, by Le Cam (1986). A gentler account of a subset of the theory appeared in the smaller book by Le Cam & Yang (1990), which will soon be reappearing in a second edition.

For too long, Le Cam's approach had an unfair reputation as dealing only with "abstract and abstruse problems" (Albers, Alexanderson & Reid 1990, page 178). But more recent work—such as the groundbreaking papers of Brown & Low (1996) and Nussbaum (1996)—has confirmed that Le Cam had long ago identified the essential mechanism at work behind the suggestive similarities between certain types of asymptotic nonparametric problems.

We can think of an experiment as a family of probability measures $\mathcal{P} = \{\mathbb{P}_\theta : \theta \in \Theta\}$ defined on a sigma-field $\mathcal{F}$ of a sample space $\Omega$. As Le Cam pointed out, the precise choice of $(\Omega, \mathcal{F})$ is somewhat irrelevant. In many respects, they are just an artifice to let us interpret the random variables of interest as measurable functions. What matters most is the ordering and vector space properties satisfied by the random variables, the properties that identify the space of all bounded, real-valued random variables on $\Omega$ as an $M$-space.

Le Cam's identification of probability measures with linear functionals on an abstract $M$-space might appear shocking, but it does make clear that there is nothing sacred about the choice of the sample space. Probability arguments that use only the $M$-space properties are valid for all choices of an underlying $\Omega$ and $\mathcal{F}$.

Strangely enough, I often find myself making analogous arguments about the irrelevance of $\Omega$ when teaching introductory probability courses. Sample spaces are very good for bookkeeping and forcing precise identification of random variables, but it is really more important to understand how those random variables relate to each other and to the probabilities. I seldom feel a need to identify explicitly an appropriate $\Omega$ when solving a problem; and even when there is an implicit $\Omega$, I would have no compunction about changing it if there were a more convenient way to represent the random variables.

Similar thoughts have occurred to me when dealing with empirical process theory. If a result is true for one sequence of random variables (or random elements of a space more complicated than the real line), and not for another sequence with the same joint distributions, but defined on a different $\Omega$, are we relying too heavily on what should be irrelevant features of the sample space? Maybe we need a stricter definition of a probability space? Such ideas are not total heresy: compare, for example, with Doob's comments on the assumption that probability measures be *perfect* in Appendix I of Gnedenko & Kolmogorov (1968), or the need for *perfect* measurable maps in the Dudley (1985) representation theorem. To paraphrase Le Cam, Why cling to an $\Omega$ that causes trouble, if there are other choices that ensure nice properties for the objects we care about?

I find myself a little too timid to take the logical next step, abandoning altogether the sample space and treating probabilities as objects that don't need the support of a sigma-field, as with Le Cam's identification of probability models with subsets of $L$-spaces in duality with $M$-spaces. It is a comfort, though, to know that there is always some $\Omega$ for which the $L$-spaces and $M$-spaces of any particular problem correspond to sets of measures and functions with an abundance of traditional regularity properties (the Kakutani representations—see pages 209–211 of Torgersen 1991).

I also see value in regarding Le Cam's more general objects as idealized measures or idealized decision procedures, particularly so if they appear only in the intermediate steps of an argument concerning their traditional counterparts. For example, a theorem asserting existence of a generalized estimator with desirable properties is a good place to start searching for a traditional estimator with the same properties, just as a mathematician who is interested only in real roots of polynomials might find it easier to prove that all the complex roots lie on the real axis, rather than establishing existence of real roots from scratch.

Let $\mathcal{Q} = \{\mathbb{Q}_\theta : \theta \in \Theta\}$ be a family of probability measures on $(\mathcal{Y}, \mathcal{A})$, with the same index set $\Theta$ as $\mathcal{P}$. Le Cam defined the distance between $\mathcal{P}$ and $\mathcal{Q}$ by means of randomizations. For example, if we have a probability kernel $K$ from $\Omega$ to $\mathcal{Y}$, then each $\mathbb{P}_\theta$ defines a probability measure $\widetilde{\mathbb{Q}}_\theta = \mathbb{P}_\theta^\omega K_\omega$ on $\mathcal{A}$, corresponding to a two-step procedure for generating an observation $y$: first generate an $\omega$ from $\mathbb{P}_\theta$, then generate a $y$ from the probability measure $K_\omega$. In more traditional notation, $\widetilde{\mathbb{Q}}_\theta(A)$ would be written $\int K_\omega(A)\,\mathbb{P}_\theta(d\omega)$, for each $A$ in $\mathcal{A}$. Again Le Cam allowed randomizations (which he called transitions) to be slightly more exotic than probability kernels, in order to achieve a most convenient compactness property for the space of all possible randomizations between two spaces—see my comments at the end of Section 3 regarding minimax bounds.

Roughly speaking, two experiments are close if each is well approximated, in the sense of total variation, by a randomization of the other. (Recall that the total variation distance between two measures $\widetilde{Q}$ and $Q$ on the same sigma-field is defined as $\|\widetilde{Q} - Q\|_1 = \sup_{|f|\le 1} |\widetilde{Q}f - Qf|$, the $f$ ranging over all measurable functions bounded in absolute value by 1. If the measures have densities $\widetilde{q}$ and $q$ with respect to a dominating measure $\mu$, then $\|\widetilde{Q} - Q\|_1$ is equal to the $L_1(\mu)$ norm of $\widetilde{q} - q$.) More formally, the distance $\delta_\Theta(\mathcal{P}, \mathcal{Q})$ is defined as the infimum over all randomizations of the quantity $\sup_{\theta \in \Theta} \|\mathbb{P}_\theta^\omega K_\omega - \mathbb{Q}_\theta\|_1$. The Le Cam distance $\Delta_\Theta(\mathcal{P}, \mathcal{Q})$ equals the minimum of $\delta_\Theta(\mathcal{P}, \mathcal{Q})$ and $\delta_\Theta(\mathcal{Q}, \mathcal{P})$. Of course, the infimum over probability kernels gives an upper bound for $\delta_\Theta(\mathcal{P}, \mathcal{Q})$, which is often all that we need.

Le Cam's method is most often applied to establish limit theory for sequences of experiments, $\mathcal{P}_n = \{\mathbb{P}_{n,\theta} : \theta \in \Theta\}$. (I will assume all $\mathbb{P}_{n,\theta}$ are defined on the same sigma-field, of the same $\Omega$, to avoid further notational complication.) If the infimum, $\delta_\Theta(\mathcal{Q}, \mathcal{P}_n)$, tends to zero as $n$ tends to infinity then I will write $\mathcal{P}_n \Subset \mathcal{Q}$, or $\mathcal{Q} \Supset \mathcal{P}_n$. If both $\mathcal{P}_n \Subset \mathcal{Q}$ and $\mathcal{P}_n \Supset \mathcal{Q}$ then the experiments $\mathcal{P}_n$ are said to converge (in Le Cam's sense) to $\mathcal{Q}$, which I will write as $\mathcal{P}_n \eqsim \mathcal{Q}$.

The randomizations involved in the definition of Le Cam's $\delta_\Theta$ fit naturally with the idea of randomized estimators (or, more generally, randomized procedures over some abstract action space). In the traditional sense, a randomized estimator of $\theta$, based on an observation $y$ from some (unknown) $\mathbb{Q}_\theta$, is just a probability kernel $\rho$ from $\mathcal{Y}$ into $\Theta$. If $K$ is a randomization from $\Omega$ to $\mathcal{Y}$, then we may define a new randomized estimator of $\theta$, based on an observation $\omega$ from $\mathbb{P}_\theta$, by averaging $\rho_y$ with respect to $K_\omega$: for each $\omega$, generate a $y$ from $K_\omega$, then generate a $t$ in $\Theta$ from $\rho_y$. I will write $K_\omega^y \rho_y$ for this new estimator, the result of composing the two randomizations.

For example, if $\mathcal{P}_n \Subset \mathcal{Q}$, with corresponding randomizations $K_n$, and if $\rho_n$ is a randomized estimator under $\mathcal{P}_n$ (that is, based on an observation $\omega$ from some $\mathbb{P}_{n,\theta}$), then the compositions $\tau_n = K_n^y \rho_{n,y}$ define a sequence of randomized estimators under $\mathcal{Q}$ (that is, based on an observation $y$ from the corresponding $\mathbb{Q}_\theta$). Any theory developed for the collection of all randomized estimators under $\mathcal{Q}$ automatically applies to $\tau_n$, and thence to $\rho_n$, within an error term derived from $\delta_\Theta(\mathcal{Q}, \mathcal{P}_n)$. See Section 3 for an example of this principle in action. Similarly, if $\mathcal{P}_n \Supset \mathcal{Q}$, then theory for randomized estimators under $\mathcal{P}_n$ has an automatic transfer to $\mathcal{Q}$.

For many important statistical results, we do not need the approximation in total variation to hold uniformly over the whole of $\Theta$. Instead, it often suffices to consider separately the finite-parameter subexperiments. That is, for each finite subset $S$ of $\Theta$, we need to find randomizations (which are allowed to depend on $S$) to make $\sup_{\theta \in S} \|\mathbb{P}_{n,\theta}^\omega K_\omega - \mathbb{Q}_\theta\|_1$ small. If $\delta_S(\mathcal{Q}, \mathcal{P}_n)$, the infimum over all such $K$, tends to zero as $n$ tends to infinity, for each finite $S$, then I will write $\mathcal{P}_n \Subset_w \mathcal{Q}$, or $\mathcal{Q} \Supset_w \mathcal{P}_n$, with the subscript $w$ standing for "weaker". Similarly, I will write $\mathcal{P}_n \eqsim_w \mathcal{Q}$ for convergence of experiments in Le Cam's weaker sense.

Convergence in Le Cam's weaker sense may be inferred from more classical results about likelihood ratios. For example, the following Lemma easily handles some of the better known results, such as the Hájek-Le Cam convolution theorems and local asymptotic minimax theorems.

<1> **Lemma.** *For $n = 1, 2, \ldots$, let $\mathcal{P}_n = \{\mathbb{P}_{n,i} : i = 1, 2, \ldots, k\}$ be a finite collection of probability measures on $(\Omega_n, \mathcal{F}_n)$, dominated by probability measures $\mathbb{P}_n$ with densities $X_{n,i} = d\mathbb{P}_{n,i}/d\mathbb{P}_n$, and let $\Omega = \{\mathbb{Q}_i : i = 1, 2, \ldots, k\}$ be probability measures on $(\mathcal{Y}, \mathcal{A})$ dominated by a probability measure $\mathbb{Q}$ with densities $Y_i = d\mathbb{Q}_i/d\mathbb{Q}$. Suppose the random vectors $X_n := (X_{n,1}, \ldots, X_{n,k})$ converge in distribution (under $\mathbb{P}_n$) to $Y := (Y_1, \ldots, Y_k)$ (under $\mathbb{Q}$). Then $\mathcal{P}_n \rightleftarrows \mathcal{Q}$.*

If $\mathbb{P}_n$ happens to be one of the $\mathbb{P}_{n,i}$, the densities $X_{n,i}$ are called *likelihood ratios*. (Likewise for $\mathcal{Q}$.) For my method of proof there is no advantage in restricting the assertion to this special case.

There are several ways to prove the Lemma. One elegant approach relates the Le Cam distance to the distance between particular measures induced on a simplex in Euclidean space by the experiments (the canonical representations of the experiments). The arguments typically proceed via a chain of results involving comparison of risk functions or Bayes risks, as in Section 2.2 of Le Cam & Yang 1990. I have no complaints about this approach, except that it can appear quite technical and forbidding to anyone wondering whether it is worth learning about Le Cam's theory. I confess that on my first few attempts at reading Le Cam (1969), the precursor to Le Cam & Yang 1990, I skipped over the relevant sections, mistakenly thinking they were of no great significance—just abstract and abstruse theory. I hope that a short, self-contained proof of the Lemma (Section 2), and an illustration of its use (Section 3), might save others from making the same mistake.

The Lemma asserts that $\mathcal{P}_n \supseteq \mathcal{Q}$ and $\mathcal{P}_n \subseteq \mathcal{Q}$. I will prove only the second assertion, which is actually the more useful. (The proof of the other assertion is almost the same.) That is, I will show that there exist probability kernels $K_{n,y}$, from $\mathcal{Y}$ to $\Omega_n$ for $n = 1, 2, \ldots$, such that $\max_i \|\mathbb{Q}_i^y K_{n,y} - \mathbb{P}_{n,i}\|_1 \to 0$ as $n \to \infty$.

My method is adapted from a more familiar part of probability theory known as the method of almost sure representation: the construction of an almost surely convergent sequence whose marginal distributions are equal to a given weakly convergent sequence of probability measures on a metric space. It is not really surprising that there should be strong similarities between Le Cam's randomizations and the coupling arguments needed to establish the representation. In both cases, one artificially constructs joint distributions as a way of deriving facts about the marginal distributions. The similarities are particularly striking for the coupling method used by Dudley (1968, 1985), who actually constructed explicitly the probability kernels needed for Lemma <1>. (Compare with the exposition in Section 9 of Pollard 1990.) The proof of the Lemma would be even shorter if I were to invoke Dudley's result directly, rather than reproducing parts of his argument to keep the proof self-contained.

## 2.  Construction of randomizations

Most calculations for Lemma <1> are carried out for a fixed $n$. To simplify notation, it helps initially to replace $\mathcal{P}_n$ by a fixed collection $\mathcal{P} = \{\mathbb{P}_i : i = 1, 2, \ldots, k\}$ of probability measures on $(\Omega, \mathcal{F})$, dominated by a fixed probability measure $\mathbb{P}$, with densities $X_i = d\mathbb{P}_i/d\mathbb{P}$.

We seek a probability kernel $K$ from $\mathcal{Y}$ to $\Omega$ to make each $\mathbb{Q}_i^y K_y$ close to the corresponding $\mathbb{P}_i$. The kernel will be chosen so that $\mathbb{Q}^y K_y = \mathbb{P}$. The kernel artificially defines a joint distribution $\mathbb{M}$ for $(y, \omega)$, with marginal distributions $\mathbb{P}$ and $\mathbb{Q}$, and conditional distribution $K_y$ for $\omega$ given $y$. The $L_1(\mathbb{M})$ distance between $X_i$ and $Y_i$, reinterpreted as random variables on $\mathcal{Y} \otimes \Omega$, will give a bound for $\|\mathbb{Q}_i^y K_y - \mathbb{P}_i\|_1$.

For the application to the proof of Lemma <1>, Assumption (iii) of the following Lemma will be checked by means of the assumed convergence in distribution of the random vectors $\{X_n\}$.

<2> **Lemma.** *Suppose there is a partition of $\mathbb{R}^k$ into disjoint Borel sets $B_0, B_1, \ldots, B_m$, and positive constants $\delta$ and $\epsilon$, for which*

  *(i) $\operatorname{diam}(B_\alpha) \leq \delta$ for each $\alpha \geq 1$*
  *(ii) $\mathbb{Q}_i\{Y \in B_0\} \leq \epsilon$ for each $i$*
  *(iii) $\mathbb{P}\{X \in B_\alpha\} \geq (1-\epsilon)\mathbb{Q}\{Y \in B_\alpha\}$ for each $\alpha$.*

*Then there exists a probability kernel $K$ from $\mathcal{Y}$ to $\Omega$ for which*
$$\max_i \|\mathbb{Q}_i^y K_y - \mathbb{P}_i\|_1 \leq 2\delta + 4\epsilon$$

*Proof.* Write $A_\alpha$ for $\{Y \in B_\alpha\}$ and $F_\alpha$ for $\{X \in B_\alpha\}$, for $\alpha = 0, 1, \ldots, m$.

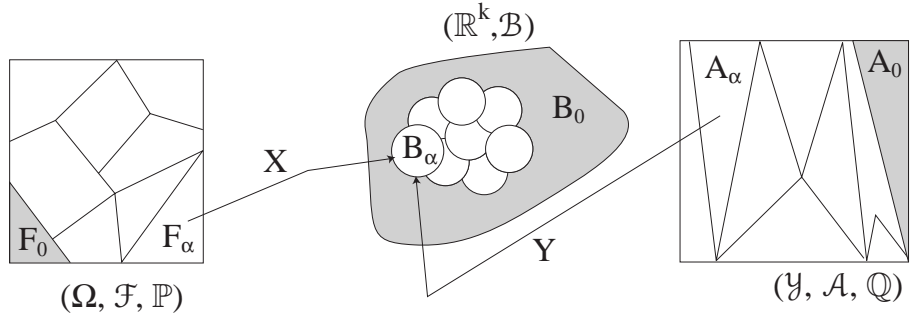

Assumption (iii) ensures that the $\mu$ defined on $\mathcal{F}$ by
$$\epsilon \mu(F) = \sum_\alpha \Big(\mathbb{P}(F_\alpha) - (1-\epsilon)\mathbb{Q}A_\alpha\Big) \mathbb{P}(F \mid F_\alpha)$$
is a probability measure. The probability kernel
$$K_y(\cdot) = \epsilon \mu(\cdot) + (1-\epsilon) \sum_\alpha \{y \in A_\alpha\} \mathbb{P}(\cdot \mid F_\alpha)$$
serves as a conditional distribution for the probability measure,

<3>
$$\mathbb{M} = \epsilon\,(\mathbb{Q} \otimes \mu) + (1-\epsilon) \sum_\alpha (\mathbb{Q}A_\alpha) \mathbb{Q}(\cdot \mid A_\alpha) \otimes \mathbb{P}(\cdot \mid F_\alpha)$$

on $\mathcal{A} \otimes \mathcal{F}$. That is,
$$\iint g(y,\omega)\,\mathbb{M}(dy, d\omega) = \int \left( \int g(y, \omega)\, K_y(d\omega) \right) \mathbb{Q}(dy),$$
at least for bounded, product measurable functions $g$; or, in the more concise linear functional notation, $\mathbb{M}^{y,\omega} g(y,\omega) = \mathbb{Q}^y K_y^\omega g(y, \omega)$. (The equality is easily checked when $g$ is the indicator function of a measurable rectangle $A \otimes F$. A generating class argument extends to more general $g$.) In particular, $\mathbb{M}$ has marginals $\mathbb{Q}$ and $\mathbb{P}$, and $\mathbb{M}Y_i(y)h(y) = \mathbb{Q}Y_i h = \mathbb{Q}_i h$ and $\mathbb{M}X_i(\omega) f(\omega) = \mathbb{P}X_i f = \mathbb{P}_i f$.

Notice that the probability measure $\mathbb{M}_0 = \sum_\alpha (\mathbb{Q}A_\alpha)\mathbb{Q}(\cdot \mid A_\alpha) \otimes \mathbb{P}(\cdot \mid F_\alpha)$ concentrates on the subset $\cup_\alpha A_\alpha \otimes F_\alpha$ of the product space. For every $(y, \omega)$ generated from $\mathbb{M}_0$, the points $X(\omega)$ and $Y(y)$ must both lie in the same $B_\alpha$. When $\alpha \geq 1$, as will happen with high $\mathbb{M}_0$ probability, the distance between $X(\omega)$ and $Y(y)$ must be smaller than $\delta$, by Assumption (i).

The Lemma asserts existence of a bound on
$$\sup_{|f| \leq 1} |\mathbb{Q}^y Y_i(y) K_y^\omega f(\omega) - \mathbb{P}^\omega X_i(\omega) f(\omega)| = \sup_{|f| \leq 1} |\mathbb{M}^{y,\omega} Y_i(y) f(\omega) - \mathbb{M}^{y,\omega} X_i(\omega) f(\omega)|$$

The right-hand side is smaller than $\mathbb{M}|Y_i - X_i|$, which we can rewrite as $2\mathbb{M}(Y_i - X_i)^+$ because $\mathbb{M}Y_i = 1 = \mathbb{M}X_i$. The contribution to $\mathbb{M}(Y_i - X_i)^+$ from $\epsilon \mathbb{Q} \otimes \mu$ is less than $\epsilon \mathbb{Q} Y_i = \epsilon$. The contribution from $(1-\epsilon)\mathbb{M}_0$ is less than
$$\mathbb{M}Y_i\{y \in A_0, \omega \in F_0\} + \mathbb{M}_0 \left(\cup_{\alpha \geq 1}\{y \in A_\alpha, \omega \in F_\alpha\}|Y_i - X_i|\right) \leq \mathbb{Q}_i A_0 + \delta.$$

☐ The assertion of the Lemma follows.

Write $Q$ for the distribution of $Y$ under $\mathbb{Q}$, a probability measure on the Borel sigma-field of $\mathbb{R}^k$. For a property like (iii) to follow from convergence in distribution, we shall need the $B_\alpha$ to be $Q$-continuity sets, that is, each boundary $\partial B_\alpha$ should have zero $Q$ measure (Pollard 1984, Section III.2). Such partitions are easy to construct from closed balls $B(y, r)$ with radius $r$ and center $y$. For each $y$, the boundaries $\partial B(y, r)$, for $0 < r < \infty$, are disjoint sets. At most countably many of those boundaries can have strictly positive $Q$ measure. In particular, for each $\delta > 0$ there is an $r_y$ with $r_y < \delta$ and $Q \partial B(y, r_y) = 0$. Some finite subcollection of such balls covers an arbitrarily large compact subset of $\mathbb{R}^k$. The partition generated by these balls then satisfies requirements (i) and (ii) of Lemma <2>.

*Partial proof of Lemma <1>.* Let $\{\delta_j\}$ and $\{\epsilon_j\}$ be sequences of positive numbers, both tending to zero. For each $j$ construct partitions $\pi_j$ over $\mathbb{R}^k$ into finite collections of $Q$-continuity sets satisfying the analog of properties (i) and (ii) of Lemma <2>, with $\delta$ replaced by $\delta_j$ and $\epsilon$ replaced by $\epsilon_j$. The convergence in distribution gives an $n_j$ for which $\mathbb{P}_n\{X_n \in B\} \geq (1 - \epsilon_j)\mathbb{Q}\{Y \in B\}$ for all $B \in \pi_j$, if $n \geq n_j$. With no loss of generality we may assume $1 = n_1 < n_2 < \ldots$. For each $n$ construct the kernel $K_{n,y}$ using the partition $\pi_j$
☐ for the $j$ such that $n_j \leq n < n_{j+1}$.

**Remarks**

(1) We don't really need $\mathbb{P}_n$ to dominate $\mathcal{P}_n$. It would suffice to have $X_{n,i}$ as the density of the part of $\mathbb{P}_{n,i}$ that is absolutely continuous with respect to $\mathbb{P}_n$, provided $\mathbb{P} X_{n,i} \to 1$ for each $i$. Such a condition is equivalent to an assumption that each $\{\mathbb{P}_{n,i} : n = 1, 2, \ldots\}$ sequence be contiguous to $\{\mathbb{P}_n\}$. Another Le Cam idea.
(2) The minimum expectation $\mathbb{M}|X_i - Y_i|$ defines the Monge-Wasserstein distance, which coincides (Dudley 1989, Section 11.8) with the bounded-Lipschitz distance between the marginal distributions. The Lemma essentially provides an alternative proof for one of the inequalities comprising Theorem 1 of Le Cam & Yang (1990, Section 2.2).

## 3. Local asymptotic minimax theorem

The theorem makes an assertion about the minimax risk for a sequence of experiments $\mathcal{P}_n = \{\mathbb{P}_{n,t} : t \in T_n\}$ for which the $T_n$ are sets that expand to some set $T$. (More formally, $\liminf T_n = T$, that is, each point of $T$ belongs to $T_n$ for all large enough $n$.) The contorted description of $T_n$ is motivated by a leading case, where $\mathbb{P}_{n,t}$ denotes the joint distribution of $n$ independent observations from some probability measure $P_{\theta_0 + A_n t}$, with $\{A_n\}$ a sequence of rescaling matrices, such as $I_d/\sqrt{n}$. If $\theta_0$ is an interior point of some subset $\Theta$ of $\mathbb{R}^d$ then the sets $T_n = \{t \in \mathbb{R}^d : \theta_0 + A_n t \in \Theta\}$ expand to the whole of $\mathbb{R}^d$, in the desired way. The experiment $\mathcal{P}_n$ captures the behavior of the $\{P_\theta : \theta \in \Theta\}$ model in shrinking neighborhoods of $\theta_0$, which explains the use of the term *local*.

In fact, the theorem has nothing to do with independence and nothing to do with the possibility that $\mathcal{P}_n$ might reflect local behavior of some other model. It matters only that $\mathcal{P}_n \in_w \mathcal{Q}$ for some experiment $\mathcal{Q} = \{\mathbb{Q}_t : t \in T\}$, a family of probability measures on some $\mathcal{Y}$. For example, $\mathcal{P}_n$ might satisfy the "local" asymptotic normality property,

$$\frac{d\mathbb{P}_{n,t}}{d\mathbb{P}_{n,0}} = (1 + \epsilon_n(t))\exp(t'Z_n - |t|^2/2) \qquad \text{for each } t \text{ in } \mathbb{R}^d,$$

where $Z_n$ is a random vector (not depending on $t$) that converges in distribution to $N(0, I_d)$ under $\mathbb{P}_{n,0}$, and $\epsilon_n(t) \to 0$ in $\mathbb{P}_{n,0}$ probability for each fixed $t$. In that case, Lemma <1> establishes $\mathcal{P}_n \ni\in_w \mathcal{Q}$ with $\mathbb{Q}_t$ denoting the $N(t, I_d)$ distribution on $\mathbb{R}^d$.

The concept of minimax risk requires a loss function $L_t(\cdot)$, defined on the space of actions. For simplicity, let me assume that the task is estimation of the parameter $t$, so that $L_t(z)$ is defined for all $t$ and $z$ in $T$. I will also assume $L_t \geq 0$. For a randomized estimator given by a probability kernel $\rho$ from $\mathcal{Y}$ to $T$, the risk is defined as

$$R(\rho, t) = \iint L_t(z) \rho_y(dz) \, \mathbb{Q}_t(dy),$$

or $\mathbb{Q}_t^y \rho_y^z L_t(z)$, when expressed in linear functional notation. The minimax risk for $\mathcal{Q}$ is defined as $R = \inf_\rho \sup_{t \in T} R(\rho, t)$, the infimum ranging over all randomized estimators from $\mathcal{Y}$ to $T$. If, for each finite constant $M$ and each finite subset $S$ of $T$, we define $R(\rho, S, M) = \max_{t \in S} \mathbb{Q}_t^y \rho_y^z (M \wedge L_t(z))$, then $R$ can be rewritten as $\inf_\rho \sup_{S,M} R(\rho, S, M)$. It is often possible (see comments below) to interchange the inf and sup, to establish a minimax equality

<4> $$R = \inf_\rho \sup_{S,M} R(\rho, S, M) = \sup_{S,M} \inf_\rho R(\rho, S, M)$$

Notice that <4> involves only the limit experiment $\mathcal{Q}$; it has nothing to do with the $\mathcal{P}_n$.

As explained by Le Cam & Yang (1990, page 84), in the context of convergence to mixed normal experiments, the following form of the theorem is a great improvement over more traditional versions of the theorem. For example, it implies that

$$\liminf_{n \to \infty} \sup_{t \in T_n} \mathbb{P}_{n,t}^\omega \tau_{n,\omega}^z L_t(z) \geq R,$$

for every sequence of randomized estimators $\{\tau_n\}$.

<5> **Theorem.** *Suppose $\mathcal{P}_n \Subset_w \mathcal{Q}$. Suppose the minimax equality <4> holds for the limit experiment $\mathcal{Q}$. Then for each $R' < R$ there exists a finite $M$ and a finite subset $S$ of the parameter set $T$ such that*

$$\inf_\tau \max_{t \in S} \mathbb{P}_{n,t}^\omega \tau_\omega^z (M \wedge L_t(z)) > R' \qquad \text{for all } n \text{ large enough.}$$

*Proof.* (Compare with Section 7.4 of Torgersen 1991.) The argument is exceedingly easy. Given constants $R' < R'' < R$, equality <4> gives a finite $M$ and a finite subset $S$ of $T$ for which

<6> $$\inf_\rho \max_{t \in S} \mathbb{Q}_t^y \rho_y^z (M \wedge L_t(z)) > R''$$

For that finite $S$, there exist randomizations $K_n$ from $\mathcal{Y}$ into $\Omega_n$ for which

$$\epsilon_n := \max_{t \in S} \| \mathbb{Q}_t^y K_{n,y} - \mathbb{P}_{n,t} \|_1 \to 0 \qquad \text{as } n \to \infty.$$

For every randomized $\tau$ under $\mathcal{P}_n$, the nonnegative function $f_t(\omega) = \tau_\omega^z (M \wedge L_t(z))$ is bounded by $M$. The function $f_t/M$ is one of those covered by the supremum that defines the total variation distance. Consequently,

$$\mathbb{P}_{n,t}^\omega \tau_\omega^z (M \wedge L_t(z)) = \mathbb{P}_{n,t}^\omega f_t(\omega) \geq \mathbb{Q}_t^y K_{n,y}^\omega f_t(\omega) - M\epsilon_n$$
$$= \mathbb{Q}_t^y K_{n,y}^\omega \tau_\omega^z (M \wedge L_t(z)) - M\epsilon_n$$

Take the maximum over $t$ in $S$. The lower bound becomes $R(\rho', S, M) - M\epsilon_n$, where $\rho' = K_{n,y}^\omega \tau_\omega$, one of the randomized estimators covered by the infimum on the left-hand side of <6>. That is, for every $\tau$, the maximum over $S$ is greater than $R'' - M\epsilon_n$, which is eventually larger than $R'$. □

As the proof showed, the theorem is not very sensitive to the definition of randomization (the $K_n$) or randomized estimator (the $\tau$ and $\rho$). It mattered only that the composition $K_{n,y}^\omega \tau_\omega$ is also a randomized estimator. With Le Cam's more general approach, with $K_n$ any transition and $\tau$ any generalized procedure, the composition is also a generalized procedure. The method of proof is identical under his approach, but with one great advantage: the minimax equality <4> comes almost for free.

Under a mild semi-continuity assumption on the loss function, there is a topology that makes the set $\mathcal{R}$ of all generalized estimators compact, and for which the map $\rho \mapsto R(\rho, S, M)$ is lower semi-continuous, for each fixed $S$ and $M$. For each $R'' < R$, the open sets $\mathcal{R}_{S,M} = \{\rho \in \mathcal{R} : R(\rho, S, M) > R''\}$ cover $\mathcal{R}$. By compactness, $\mathcal{R}$ has a finite subcover $\cup_{i=1,\ldots,m} \mathcal{R}_{S_i, M_i}$. With $S = \cup_i S_i$ and $M = \max_i M_i$, we have $R(\rho, S, M) > R''$ for every $\rho$ in $\mathcal{R}$, as required for <4>.

In effect, by enlarging the class of objects regarded as randomizations or randomized estimators, Le Cam cleverly removed the *minimax* from the hypotheses of the local asymptotic minimax theorem. With the right definitions, the theorem is an immediate consequence of the lower semi-continuity of risk functions—which is essentially what Le Cam (1986, Section 7.4) said.

**Acknowledgement.**

In correspondence during the summer of 1990, I discussed some of the ideas appearing in this paper with Lucien Le Cam. Even when he disagreed with my suggestions, or when I was "discovering" results that were already in his 1986 book or earlier work, he was invariably gracious, encouraging and helpful.

REFERENCES

Albers, D., Alexanderson, J. & Reid, C. (1990), *More Mathematical People*, Harcourt, Brace, Jovanovich, Boston.

Brown, L. D. & Low, M. G. (1996), 'Asymptotic equivalence of nonparametric regression and white noise', *Annals of Statistics* **24**, 2384–2398.

Dudley, R. M. (1968), 'Distances of probability measures and random variables', *Annals of Mathematical Statistics* **39**, 1563–1572.

Dudley, R. M. (1985), 'An extended Wichura theorem, definitions of Donsker classes, and weighted empirical distributions', *Springer Lecture Notes in Mathematics* **1153**, 141–178. Springer, New York.

Dudley, R. M. (1989), *Real Analysis and Probability*, Wadsworth, Belmont, Calif.

Gnedenko, B. V. & Kolmogorov, A. N. (1968), *Limit Theorems for Sums of Independent Random Variables*, Addison-wesley.

Le Cam, L. (1969), *Théorie Asymptotique de la Décision Statistique*, Les Presses de l' Université de Montréal.

Le Cam, L. (1986), *Asymptotic Methods in Statistical Decision Theory*, Springer-Verlag, New York.

Le Cam, L. & Yang, G. L. (1990), *Asymptotics in Statistics: Some Basic Concepts*, Springer-Verlag.

Nussbaum, M. (1996), 'Asymptotic equivalence of density estimation and gaussian white noise', *Annals of Statistics* **24**, 2399–2430.

Pollard, D. (1984), *Convergence of Stochastic Processes*, Springer, New York.

Pollard, D. (1990), *Empirical Processes: Theory and Applications*, Vol. 2 of *NSF-CBMS Regional Conference Series in Probability and Statistics*, Institute of Mathematical Statistics, Hayward, CA.

Torgersen, E. (1991), *Comparison of Statistical Experiments*, Cambridge University Press.